\documentclass[12pt,english,a4paper,oneside]{amsart}

\usepackage[cp1251]{inputenc}
\usepackage[english]{babel}
\usepackage{amsmath,amsthm,amssymb}

\usepackage[
a4paper,
left=22mm,
right=18mm,
top=20mm,
bottom=27mm
]{geometry}

\usepackage{xcolor}
\usepackage[colorlinks]{hyperref}

\newtheorem{theorem}{Theorem}[section]
\newtheorem{lemma}[theorem]{Lemma}
\newtheorem{corollary}[theorem]{Corollary}

\numberwithin{equation}{section}
\theoremstyle{definition}

\title{Empirical approximation of $L_p$ norms for  $p\in(0,2)$}

\author{Feng Dai}
\address{F.~Dai, Department of Mathematical and Statistical Sciences, University of Alberta, Edmonton, Alberta T6G 2G1, Canada}
\email{fdai@ualberta.ca}

\author{Egor Kosov}
\address{E.~Kosov, Centre de Recerca Matem\`atica, Campus de Bellaterra, Edifici~C 08193
Bellaterra (Barcelona), Spain.}
\email{kosoved09@gmail.com}

\thanks{The first author is supported by the
NSERC of Canada Discovery Grant (Grant No. RGPIN 2026-05667).
The second author is supported by 
by the AEI grants 
RYC2023-043616-I and
PID2025-169712NA-I00 funded by MICIU/AEI/10.13039/501100011033,
and by the Spanish State Research Agency, through the Severo Ochoa and Mar\'ia de Maeztu Program for Centers and
Units of Excellence in R\&D (CEX2020-001084-M).
He thanks CERCA Programme (Generalitat de Catalunya) for institutional support.
}

\begin{document}
\maketitle

\begin{abstract}
We introduce a new approach to
proving sharp large deviation inequalities for
suprema of random power sums of the form
\[
\sum_{j=1}^m\zeta_j w_j
|\langle{\bf v}_j,y\rangle|^p,\quad p\in(0,2),
\]
where $\zeta_1,\ldots,\zeta_m$ are independent subgaussian random variables.

We apply this method to show that,
for any $N$-dimensional subspace $F_N$
of bounded functions on a probability space $(\Omega, \mu)$
such that the natural embedding
$(F_N,\|\cdot\|_{L_2(\mu)})\hookrightarrow L_\infty(\mu)$
has operator norm $O(\sqrt N)$,
$m=O(N\log N)$ independent samples $\xi^1,\ldots,\xi^m$ drawn from $\mu$
suffice to ensure the simultaneous sampling discretization
\[
(1-\varepsilon)\|f\|_{L_p(\mu)}^p
\le\frac1m\sum_{j=1}^m|f(\xi^j)|^p
\le(1+\varepsilon)\|f\|_{L_p(\mu)}^p
\quad\forall f\in F_N\quad\forall p\in[1,2]
\]
with high probability.
The order $N\log N$ cannot be improved in general
for independent random sampling points.

In particular, this result resolves Talagrand's question on random embeddings
of subspaces of $L_p$, $p\in(1,2)$, into $\ell_p^m$.
Namely, every $N$-dimensional subspace of $L_p$,
$p\in(1,2)$, admits an embedding into $\ell_p^m$
obtained by random sampling of $m$ points after a Lewis change of density,
with distortion at most
$1+\varepsilon$ and
$m\le C\varepsilon^{-2}N\log N$.
\end{abstract}

\section{Introduction}

\subsection{Main results}

Throughout this paper, a real-valued random variable $\zeta$ is called  $\sigma$-subgaussian for some constant $\sigma>0$ if it is centered and satisfies
\[
\mathbb E e^{t\zeta}\le e^{\sigma^2 t^2/2}
\quad \text{for every}\quad  t\in\mathbb R.
\]

Our first main result is the following uniform moment estimate
for sums of independent subgaussian random variables indexed
by powers of vector coordinates.

\begin{theorem}\label{th-uniform-signed-sum}
Let $0<p_0\le1$ and let
$E\subset\mathbb R^m$ be an $N$-dimensional subspace satisfying
\[
\|x\|_\infty\le\sqrt\delta\,\|x\|_2
\quad  \text{for every}\quad    x\in E.
\]
If $\zeta_1,\ldots,\zeta_m$ are independent $1$-subgaussian random
variables, then, for every $r\ge1$,
\[
\biggl[
\mathbb E
\biggl(
\sup_{p_0\le p\le2}\sup_{x\in E\cap B_p^m}
\Bigl|\sum_{j=1}^m\zeta_j|x_j|^p\Bigr|
\biggr)^r
\biggr]^{1/r}
\le
8\sqrt{\delta\bigl(p_0^{-1}\log (2N)+r\bigr)},
\]
where $B_p^m=\{x\in\mathbb R^m\colon \sum_{j=1}^m|x_j|^p\le1\}$.
\end{theorem}

For fixed $p=2$, the estimate is   due to Rudelson \cite{Rud,Rud99}, while
the corresponding estimate for fixed $p=1$ goes back
to Talagrand \cite{Tal90}.
For $1<p<2$, Talagrand \cite{Talagrand1995} obtained
a weaker estimate with an additional factor involving
$\log(\delta^{-1})$ and $\log\log m$.
For $0<p<1$, analogous estimates, with similar logarithmic losses, were
proved by Zvavitch \cite{Zva07}.
Theorem \ref{th-uniform-signed-sum} removes these logarithmic losses and yields  a bound
that holds simultaneously for all $p\in[p_0,2]$.
For each fixed $p_0>0$, the dependence on $\delta$, $N$, and $r$
is sharp, as can be seen by taking independent standard Gaussian
coefficients and letting $E$ consist of vectors that are constant
on each of $N$ disjoint blocks of equal size.

An important application of Theorem~\ref{th-uniform-signed-sum}
concerns the classical problem of almost-isometric embeddings of
$N$-dimensional subspaces of $L_p[0,1]$ into $\ell_p^m$.
This problem goes back to Johnson and Schechtman
\cite{JS82,Schechtman1987} (see \cite{JS} for a detailed discussion).
Given $N$, $p$, and $\varepsilon>0$, one seeks the smallest
dimension $m=m(N,p,\varepsilon)$ such that every such subspace
admits a linear embedding into $\ell_p^m$ with distortion at most
$1+\varepsilon$. For even integers $p$, Schechtman~\cite{Schechtman2011}
derived the optimal order in $N$ from the Batson--Spielman--Srivastava sparsification
theorem~\cite{BSS2012}. For $p=1$, the optimal order  was recently obtained
by Reis and Rothvoss~\cite{ReisRothvoss2026}.
For general $p>2$, Bourgain, Lindenstrauss, and Milman
\cite{BLM1989} proved that one can take
\[
m\le C(p,\varepsilon)N^{p/2}\log N.
\]

For $1<p<2$, Talagrand~\cite{Talagrand1995} established the bound
\[
m\le C(p,\varepsilon)N\log N(\log\log N)^2.
\]
His proof relies on a technically involved generic chaining
argument, and he later described the theorem as
``one of the most difficult (and pretty) results'' of his book
\cite[p.~573]{Talagrand2014}. The constant in the estimate deteriorates
as $p\uparrow2$ (see \cite[p.~572]{Talagrand2014}). 
Talagrand explicitly identified $N\log N$ as his original target
for the random-choice approach he used and described the additional
iterated-logarithm factor $(\log\log N)^2$ as ``truly parasitic''
\cite[pp.~311--312]{Talagrand1995}. The question of removing this factor was subsequently
discussed in
\cite[discussion following Theorem~16.8.1]{Talagrand2014}
and \cite[discussion following Lemma~7.3]{CohenPeng2015}.

In this paper, we resolve this question on random embeddings, which was posed more
than thirty years ago.
Namely, we show 
that for every
$0<p<2$ and $0<\varepsilon<1$, each $N$-dimensional subspace
of $L_p$ embeds linearly into $\ell_p^m$ with distortion at most
$1+\varepsilon$ and
\[
m\le Cp^{-3}\varepsilon^{-2}N\log N,
\]
where $C>0$ is an absolute constant
(see Corollary~\ref{cor-Lp-embedding} below). In particular, the estimate avoids the deterioration
of Talagrand's bound as $p\uparrow2$.
For $0<p<1$, our result also  removes the factor
$(\log\log N)^2$ from the bound on the embedding dimension in
Zvavitch's theorem \cite[Theorem~1.1]{Zva07}.

The embedding is constructed from independent,
identically distributed samples drawn after Lewis's change of density.
This reduces the construction to the following problem of sampling
discretization of $L_p$ norms.
For an $N$-dimensional space $F_N$ of bounded measurable functions
on a probability space $(\Omega,\mathcal F,\mu)$
and $\varepsilon\in(0,1)$, we seek the smallest number $m$
for which there exist points $\xi^1,\ldots,\xi^m\in\Omega$ such that
\begin{equation}\label{eq-mz}
(1-\varepsilon)\|f\|_{L_p(\mu)}^p
\le
\frac{1}{m}\sum_{j=1}^m |f(\xi^j)|^p
\le
(1+\varepsilon)\|f\|_{L_p(\mu)}^p
\quad\forall f\in F_N.
\end{equation}
To define point evaluations, we fix bounded measurable
representatives of a basis of $F_N$ and use their linear combinations.

A natural approach is to choose the points independently
with distribution $\mu$ and ask whether the discretization
inequality \eqref{eq-mz} holds with high probability.
For $p=2$, under the natural normalization assumption
\[
\|f\|_\infty:=\sup_{x\in\Omega}|f(x)|
\le\sqrt{KN}\,\|f\|_{L_2(\mu)}
\quad\forall f\in F_N,
\]
such discretization results follow from the matrix Chernoff
inequality \cite[Corollary~5.2]{Tropp2012}
whenever $m\ge CKN\varepsilon^{-2}\log N$
(see also \cite[Theorem~3.2]{DKT26}).
Bounds on the expected uniform error follow from earlier work
of Rudelson \cite{Rud99}.

For $p\in(1,2)$, a modification of Talagrand's
argument \cite{Talagrand1995} was used in \cite{DKT} to prove,
under the same normalization assumption, the existence of
a discretization set with cardinality of order
$C(p,K,\varepsilon)N\log N(\log\log N)^2$.
Theorem \ref{th-uniform-signed-sum} allows us to remove the additional factor $(\log\log N)^2$.
For every fixed $p\in[1,2]$, $K\ge1$, and $\varepsilon\in(0,1)$,
our result yields the optimal sample-size order $N\log N$
for independent random sampling.
Moreover, we obtain a single sampling set that simultaneously
discretizes all the $L_p$ norms with $p\in[p_0,2]$,
for each fixed $p_0>0$.

\begin{theorem}\label{th-uniform-sampling}
Let $p_0\in(0, 1]$,
$C:=2^{21}$, $c:=2^{-20}$, $K\ge1$, and let
$F_N$ be a real or complex $N$-dimensional space of
bounded measurable functions on a probability space $(\Omega,\mu)$.
Assume that
\[
\|f\|_\infty\le\sqrt {KN}\,\|f\|_{L_2(\mu)}
\quad\forall f\in F_N.
\]
Then, for every $\varepsilon\in(0,1)$ and every $m\in\mathbb N$ satisfying
\[
m\ge Cp_0^{-1} \varepsilon^{-2}KN\log(2N),
\]
independent points $\xi^1,\ldots,\xi^m$ with common distribution $\mu$
satisfy
\[
(1-\varepsilon)\|f\|_{L_p(\mu)}^p
\le\frac1m\sum_{j=1}^m|f(\xi^j)|^p
\le(1+\varepsilon)\|f\|_{L_p(\mu)}^p
\quad\forall f\in F_N\quad \forall p\in[p_0, 2]
\]
with probability at least
$1-\exp\bigl(-c\frac{m\varepsilon^2}{KN}\bigr)$.
\end{theorem}

The optimality of the sample-size order  $N\log N$ for independent random
sampling follows by considering the space of all functions on
an $N$-point set equipped with the uniform probability measure
(see, for instance, \cite[Example~4.1]{DKM}).
For further discussion of sampling discretization, see
the surveys \cite{DKT26, KKLT}.
Related results for $p>2$ were obtained in \cite{DKM}.

We finally note that the same sampling estimate improves the sample-size
guarantees for the Lewis-weight row-sampling method of
Cohen and Peng \cite{CohenPeng2015}.
For fixed $p\in(1,2)$, fixed accuracy $\varepsilon\in(0,1)$,
and fixed success probability, it removes the factor
$(\log\log N)^2$ from their bound, so that $O(N\log N)$
sampled rows suffice, where $N$ is the rank of the input matrix.

\subsection{Proof strategy}
To prove Theorem~\ref{th-uniform-signed-sum}, we introduce
a new polynomial lifting approach to random power sums
of the form
\begin{equation}\label{eq-process}
\sum_{j=1}^m \zeta_j w_j
|\langle{\bf v}_j,y\rangle|^p = 
\sum_{j=1}^m \zeta_j w_j
|\langle{\bf v}_j{\bf v}_j^* y,y\rangle|^{p/2},\quad p\in(0, 2].
\end{equation}
This approach differs from Talagrand's generic chaining
and from other methods based on metric entropy estimates
and relies on reducing the problem to estimating random matrix sums.
The key step is to represent the fractional powers
$|\langle{\bf v}_j{\bf v}_j^* y,y\rangle|^{p/2}$
as positively weighted series of quadratic forms on enlarged spaces
of homogeneous polynomials. This expresses the nonlinearly indexed process \eqref{eq-process}
in terms of random sums of self-adjoint operators, allowing us
to use the machinery of random matrix theory.
In particular, we use the Lieb--Tropp matrix Laplace
inequality \cite[Corollary 3.4.2 and Lemma~3.5.1]{Tropp2015} to obtain probability
estimates for these random sums.
The central difficulty is that the dimensions of the lifted
spaces grow rapidly with the polynomial degree, so the
dimension factors in standard matrix concentration bounds
are too large to give the desired estimates.
Thus, we combine this lifting procedure with new trace estimates
for the lifted operators associated with the corresponding
deterministic power sums.
These estimates exploit the structure of the operators to
control the traces of their exponentials without incurring
the full dimension factor of the lifted space.

Beyond the applications discussed here, we expect the proposed
approach to extend to more general random sums of the form
\[
\sum_{j=1}^m \langle A_jy,y\rangle^{p/2},
\]
where $A_1,\ldots,A_m$ are independent random positive
semidefinite matrices.
Such an extension would generalize matrix Chernoff bounds in Tropp's form
\cite[Corollary~5.2]{Tropp2012} to fractional powers of
quadratic forms and yield improved bounds for the
sparsification of sums of powers of Euclidean norms,
a class of problems considered in
\cite{JLLS2023, ManojOvsiankin2025}.

\subsection{Notation}

For $0<p\le\infty$, $B_p^m$ denotes the closed unit ball
of $\ell_p^m$, with its usual quasi-norm when $p<1$.

For a linear operator $A$ on a finite-dimensional Hilbert space,
$A^*$ denotes its adjoint, $\operatorname{tr}A$ its trace,
and $\|A\|_{\mathrm{op}}$ and $\|A\|_{\mathrm{HS}}$
its operator and Hilbert--Schmidt norms, respectively.
Thus,
\[
\|A\|_{\mathrm{HS}}^2=\operatorname{tr}(AA^*).
\]
We denote the identity operator on a space $H$ by $I_H$,
and write $I_N$ when $H=\mathbb R^N$.

For self-adjoint operators $A$ and $B$, we write $A\preceq B$
if
\[
\langle Ax,x\rangle\le\langle Bx,x\rangle
\quad\text{for every }x.
\]
Subscripts on $\mathbb E$ indicate the random variables
with respect to which the expectation is taken.


\subsection{Organization of the paper}

In Section~\ref{sect-pol}, we introduce the main technique
underlying the proofs of our results.
In Section~\ref{sect-power-sums}, we prove
Theorem~\ref{th-uniform-signed-sum}.
In Section~\ref{sect-sampling}, we deduce
Theorem~\ref{th-uniform-sampling}, from which we obtain
the stated results on embeddings of subspaces of $L_p$
in Section~\ref{sect-embedding}.

\section{Lifting to the space of homogeneous polynomials}
\label{sect-pol}

\subsection{The main construction}

We first explain the main idea behind the proof of
Theorem~\ref{th-uniform-signed-sum}
by considering a simplified version without the supremum over $p$.
For fixed $p\in(0,2]$,
our goal is to estimate the moments of
\[
\sup_{x\in E\cap B_p^m}
\Bigl|\sum_{j=1}^m\zeta_j|x_j|^p\Bigr|.
\]
By choosing a basis of $E$, we reduce the problem
to estimating the moments of
\[
\sup_{\|y\|_{\ell^p_w}\le 1}
\Bigl|\sum_{j=1}^m\zeta_jw_j
|\langle{\bf v}_j,y\rangle|^p\Bigr|,
\]
where
${\bf v}_1,\ldots,{\bf v}_m\in\mathbb R^N$ are unit vectors,
$w_1,\ldots,w_m\ge0$ are weights that may depend on $p$,
and where
\[
\|y\|_{\ell^p_w}:=
\biggl(\sum_{j=1}^m w_j|\langle{\bf v}_j,y\rangle|^p\biggr)^{1/p}
\quad y\in\mathbb R^N.
\]
Let $\alpha = p/2\in (0, 1]$.
By homogeneity
\[
\sup_{\|y\|_{\ell^p_w}\le 1}
\Bigl|\sum_{j=1}^m\zeta_jw_j
|\langle{\bf v}_j,y\rangle|^p\Bigr|
=\sup_{\|y\|_2=1}
\frac{
\bigl|\sum_{j=1}^m\zeta_j w_j
|\langle{\bf v}_j,y\rangle|^{2\alpha}\bigr|
}{
\bigl|\sum_{j=1}^m w_j
|\langle{\bf v}_j,y\rangle|^{2\alpha}\bigr|
}.
\]
When $\alpha=1$, the sum depends linearly on the rank-one
operators ${\bf v}_j{\bf v}_j^*$, and we can write
\[
\sum_{j=1}^m\zeta_j w_j\langle{\bf v}_j,y\rangle^2
=
\biggl\langle
\Bigl(\sum_{j=1}^m\zeta_j w_j{\bf v}_j{\bf v}_j^*\Bigr)y,y
\biggr\rangle.
\]
This allows us to work with random sums of operators.
The difficulty in the case $\alpha<1$ is the lack of this
linearity.
The key observation is that the powers
$|\langle{\bf v}_j,y\rangle|^{2\alpha}$
can be represented as positively weighted series of {\it quadratic forms} on spaces of homogeneous polynomials.

Indeed, for every $\alpha\in(0,1]$ and every $t\in[0,1]$,
the Taylor expansion gives
\[
(1-t)^\alpha
=
1-\sum_{k=1}^\infty b_{\alpha,k}t^k,
\quad
b_{\alpha,k}=(-1)^{k+1}\binom\alpha k\ge0.
\]
In particular,
\begin{equation}\label{eq-coeff-equality}
\sum_{k=1}^\infty b_{\alpha,k}=1
\end{equation}
and, substituting $t=1-s$, we obtain
\[
s^\alpha
=
\sum_{k=1}^\infty b_{\alpha,k}[1-(1-s)^k]
\quad\forall s\in[0,1].
\]
Hence, for every unit vector $y\in\mathbb R^N$,
\[
|\langle{\bf v}_j,y\rangle|^{2\alpha}
=
\sum_{k=1}^\infty b_{\alpha,k}
\bigl[1-(1-\langle {\bf v}_j,y\rangle^2)^k\bigr].
\]
We now represent each expression
$1-(1-\langle {\bf v}, y\rangle^2)^k$
as a quadratic form associated with an orthogonal projection
on the space $\mathcal H_k^N$ of all homogeneous polynomials of degree $k$
in $N$ variables
$z=(z_1,\ldots,z_N)\in\mathbb R^N$.
We equip this space with the Bombieri inner product
\[
\biggl\langle
\sum_{|\beta|=k}a_\beta z^\beta,
\sum_{|\beta|=k}b_\beta z^\beta
\biggr\rangle_{\mathcal H_k^N}
:=
\sum_{|\beta|=k}
\frac{a_\beta b_\beta}{\binom{k}{\beta}},
\quad
\binom{k}{\beta}:=\frac{k!}{\beta_1!\ldots\beta_N!},
\quad
\beta=(\beta_1,\ldots,\beta_N).
\]
In particular, the normalized monomials
$\bigl\{\binom{k}{\beta}^{1/2}z^\beta\colon|\beta|=k\bigr\}$
form an orthonormal basis with respect to this inner product.

For $y\in\mathbb R^N$, let
\[
\psi_{k,y}(z):=\langle y,z\rangle^k.
\]
For a linear operator $B\colon\mathbb R^N\to\mathbb R^N$, let
\[
\Gamma_k(B)\colon\mathcal H_k^N\to\mathcal H_k^N,
\quad
(\Gamma_k(B)f)(z):=f(B^*z).
\]
For a unit vector ${\bf v}\in\mathbb R^N$, we define the operator
\[
T_{k,{\bf v}}
:=
I_{\mathcal H_k^N}-\Gamma_k(P_{{\bf v}^\perp}),
\]
where, for a subspace $H\subset\mathbb R^N$,
$P_H$ denotes the orthogonal projection onto $H$.
The key property of these operators, which we verify below,
is the identity
\begin{equation}\label{eq-key-identity}
\langle\psi_{k,y},T_{k,{\bf v}}\psi_{k,y}\rangle_{\mathcal H_k^N}
=
1-(1-\langle {\bf v}, y\rangle^2)^k
\quad
\forall y,{\bf v}\in\mathbb R^N
\text{ with }\|y\|_2=\|{\bf v}\|_2=1.
\end{equation}
Thus,
\[
|\langle{\bf v}_j,y\rangle|^{2\alpha}
=
\sum_{k=1}^\infty b_{\alpha,k}
\langle\psi_{k,y},T_{k,{\bf v}_j}\psi_{k,y}\rangle_{\mathcal H_k^N}
,
\]
and, for arbitrary weights $w_1, \ldots, w_m$,
we obtain
\begin{equation}\label{eq-repres}
\sum_{j=1}^m w_j|\langle {\bf v}_j, y \rangle|^{2\alpha}
\!=\!\sum_{k=1}^\infty b_{\alpha, k}\!
\sum_{j=1}^m w_j
\langle\psi_{k,y},T_{k, {\bf v}_j}\psi_{k,y}\rangle_{\mathcal H_k^N}
\!=\! \sum_{k=1}^\infty b_{\alpha, k}\Bigl\langle
\psi_{k,y}, \Bigl(\sum_{j=1}^m w_jT_{k, {\bf v}_j}\Bigr)\psi_{k,y}
\Bigr\rangle_{\mathcal H_k^N}
\end{equation}
for every unit vector $y\in\mathbb R^N$.
Similarly, for random sums,
\begin{equation}\label{eq-repres-random}
\sum_{j=1}^m\zeta_j w_j
|\langle{\bf v}_j,y\rangle|^{2\alpha}
=
\sum_{k=1}^\infty b_{\alpha,k}
\Bigl\langle
\psi_{k,y},
\Bigl(\sum_{j=1}^m\zeta_j w_jT_{k,{\bf v}_j}\Bigr)\psi_{k,y}
\Bigr\rangle_{\mathcal H_k^N}
\quad
\forall y\in\mathbb R^N\text{ with }\|y\|_2=1.
\end{equation}

The representation \eqref{eq-repres-random} expresses the random
power sum through quadratic forms of random sums of self-adjoint
operators.
This allows us to use the machinery of random matrix
theory to obtain probability estimates.
In particular, we will use the Lieb--Tropp matrix Laplace
inequality \cite[Corollary~3.4.2]{Tropp2015}, which asserts that
\begin{equation}\label{eq-LT}
\mathbb E\operatorname{tr}\exp\bigl(H+A\bigr)
\le 	
\operatorname{tr}\exp\bigl(
H+\log\mathbb E e^{A}
\bigr),
\end{equation}	
whenever $H$ is a deterministic self-adjoint matrix and $A$ is a random self-adjoint matrix
of the same size with $e^{A}$ integrable.

We use the Lieb--Tropp matrix Laplace
inequality to control the random operators
$\sum_{j=1}^m \zeta_jw_jT_{k,{\bf v}_j}$
in terms of the deterministic operators
\[
S_k:=\sum_{j=1}^m w_jT_{k,{\bf v}_j}
\]
with high probability. The resulting probability estimates are then reduced
to bounds on $\operatorname{tr}e^{-tS_k}$.

\subsection{Proof of the identity \eqref{eq-key-identity}}

We now verify several simple properties of $\Gamma_k$ and
$T_{k,{\bf v}}$ and establish the key identity \eqref{eq-key-identity}.

\begin{lemma}\label{lem-aux}
The following properties of the mapping $\Gamma_k$ hold:

$(i)\, \Gamma_k(B_1 B_2) = \Gamma_k(B_1)\Gamma_k(B_2)$.

$(ii)\, (\Gamma_k(B))^* = \Gamma_k(B^*)$.

$(iii)$ 
For every subspace $H\subset\mathbb R^N$,
the operator $\Gamma_k(P_H)$, where $P_H$ denotes the
orthogonal projection onto $H$, is also an orthogonal projection.

$(iv)$ The operator
$T_{k, {\bf v}} := I_{\mathcal{H}_k^N} - \Gamma_k(P_{{\bf v}^\perp})$ is an orthogonal projection and
\[
\langle\psi_{k,y},T_{k, {\bf v}}\psi_{k,y}\rangle_{\mathcal H_k^N} = 1-(1-\langle {\bf v}, y\rangle^2)^k
\quad \forall y, {\bf v}\in \mathbb{R}^N\text{ with } \|y\|_2=\|{\bf v}\|_2=1.
\]
\end{lemma}

\begin{proof}
The property $(i)$ follows directly from the definition.	

Expanding the power, we obtain
\[
\psi_{k, y}(z)
=\biggl(\sum_{j=1}^N y_jz_j\biggr)^k
=\sum_{|\beta| = k} \binom{k}{\beta}y^\beta z^\beta.
\]	
Therefore,
\begin{equation}\label{eq-psi-equality}
\langle\psi_{k, y},\psi_{k, y'}\rangle_{\mathcal H_k^N}
=\sum_{|\beta| = k} \binom{k}{\beta}y^\beta (y')^\beta = \langle y,y'\rangle^k.
\end{equation}
In particular, $\|\psi_{k,y}\|_{\mathcal H_k^N}=1$ whenever $\|y\|_2=1$.
Moreover,
\[
f(y) = \langle\psi_{k,y}, f\rangle_{\mathcal H_k^N}
\quad \forall f\in \mathcal H_k^N.
\]	
Thus,
\begin{equation}\label{eq-span}
\operatorname{span}\{\psi_{k,y}\colon y\in \mathbb{R}^N\} = \mathcal H_k^N,
\end{equation}
since only the zero polynomial can be orthogonal
to every $\psi_{k,y}$.

Since
\[
\Gamma_k(B)\psi_{k,y} = \psi_{k,By},
\]
by \eqref{eq-psi-equality},
we obtain
\[
\bigl\langle (\Gamma_k(B))^*\psi_{k,y},\psi_{k,y'}\bigr\rangle_{\mathcal H_k^N}
=
\bigl\langle \psi_{k,y}, \psi_{k, By'}\bigr\rangle_{\mathcal H_k^N}
=\langle y, By'\rangle^k
=\langle B^*y, y'\rangle^k
=
\bigl\langle \Gamma_k(B^*)\psi_{k,y}, \psi_{k, y'}\bigr\rangle_{\mathcal H_k^N}.
\]
The property $(ii)$ now follows from \eqref{eq-span}.

By $(i)$ and $(ii)$,
\[
\Gamma_k(P_H)^2
=\Gamma_k(P_H^2)
=\Gamma_k(P_H)
\quad\text{and}\quad
\Gamma_k(P_H)^*
=\Gamma_k(P_H^*)
=\Gamma_k(P_H),
\]
which proves $(iii)$.

Finally,
$T_{k, {\bf v}}$ is an orthogonal projection
by $(iii)$,
and, by \eqref{eq-psi-equality},
\[
\langle\psi_{k,y},T_{k, {\bf v}}\psi_{k,y}\rangle_{\mathcal H_k^N}
=1-\langle y, P_{{\bf v}^\perp}y\rangle^k
=1-\bigl(1-\langle{\bf v},y\rangle^2\bigr)^k,
\]
which proves $(iv)$.
\end{proof}

\subsection{A trace estimate for isotropic frames}

As explained above,
the key technical ingredient in
controlling the random operators in \eqref{eq-repres-random}
with high probability
using the Lieb--Tropp matrix Laplace
inequality
is the following trace estimate.

\begin{lemma}\label{lem-heat-trace}
Let
${\bf v}_1,\ldots,{\bf v}_m\in\mathbb R^N$ and
$w_1,\ldots,w_m\ge0$ satisfy
\[
\|{\bf v}_j\|_2=1
\quad\text{and}\quad
\sum_{j=1}^m w_j{\bf v}_j{\bf v}_j^*=I_N.
\]
Let $k\in\mathbb{N}$, and let
\[
S_k:=\sum_{j=1}^m w_jT_{k, {\bf v}_j},
\quad T_{k, {\bf v}_j} := I_{\mathcal{H}_k^N} - \Gamma_k(P_{{\bf v}_j^\perp}).
\]
Then
\[
\operatorname{tr}e^{-tS_k}
\le
Ne^{-t/2} e^{Nke^{-t/2}}\quad \forall t>0.
\]
\end{lemma}

\begin{proof}
By the assumptions, $\sum_{j=1}^mw_j=N$.	
Let $J_1, J_2, \ldots$
be independent random variables with the distribution
\[
\mathbb P(J_n=j)=\frac{w_j}{N}\quad j\in\{1, \ldots, m\},
\]
and, for a
fixed $t\ge0$, let $\xi$ be a Poisson random variable
with mean $Nt/2$, independent of the sequence $\{J_n\}_{n\ge 1}$.
Set
\[
B_0=C_0=I_N,
\quad
B_n:=P_{{\bf v}_{J_n}^\perp}\cdot\ldots\cdot P_{{\bf v}_{J_1}^\perp},
\quad
C_n:=B_nB_n^*.
\]
By Lemma \ref{lem-aux}(i) and independence,
we have
\[
\mathbb E\Gamma_k(B_n)
=\mathbb E\bigl[\Gamma_k(P_{{\bf v}_{J_n}^\perp})\cdot\ldots\cdot \Gamma_k(P_{{\bf v}_{J_1}^\perp})\bigr]
=\bigl(\mathbb E\Gamma_k(P_{{\bf v}_{J_1}^\perp})\bigr)^n
=\Bigl(I_{\mathcal H_k^N}-\frac{S_k}{N}\Bigr)^n,
\]
and therefore
\[
\mathbb E\Gamma_k(B_\xi)
=e^{-Nt/2}\sum_{n=0}^\infty
\frac{(Nt/2)^n}{n!}
\Bigl(I_{\mathcal H_k^N}-\frac{S_k}{N}\Bigr)^n
=e^{-tS_k/2}.
\]
Since the operator $e^{-tS_k/2}$ is self-adjoint,
\begin{align}\label{eq-1-trace-est}
\operatorname{tr}e^{-tS_k}
&=\|e^{-tS_k/2}\|_{\mathrm{HS}}^2
=\|\mathbb E\Gamma_k(B_\xi)\|_{\mathrm{HS}}^2
\le\mathbb E\|\Gamma_k(B_\xi)\|_{\mathrm{HS}}^2
\\
&=
e^{-Nt/2}\sum_{n=0}^\infty
\frac{(Nt/2)^n}{n!}\mathbb E\|\Gamma_k(B_n)\|_{\mathrm{HS}}^2.\nonumber
\end{align}
We now estimate $\|\Gamma_k(B_n)\|_{\mathrm{HS}}^2$.

For a positive semidefinite linear operator $C\colon \mathbb{R}^N\to \mathbb{R}^N$ with eigenvalues $\lambda_1, \ldots, \lambda_N$,
let
\[
e_\ell(C):=
\sum_{\substack{J\subset\{1,\ldots,N\}\\|J|=\ell}}
\prod_{j\in J}\lambda_j.
\]
In particular, $e_0(C)=1$ and
$e_\ell(I_N)=\binom N\ell$.
The eigenvalues of $\Gamma_k(C)$ are the monomials
$\lambda_1^{\beta_1}\cdots \lambda_N^{\beta_N}$ with $|\beta| = k$, and therefore
\[
\operatorname{tr}(\Gamma_k(C))
=\sum_{|\beta|=k}\lambda_1^{\beta_1}\ldots
\lambda_N^{\beta_N}.
\]
When $\|C\|_{\rm op}\le1$, grouping the terms according to their supports, we obtain
\begin{align}\label{eq-gamma-hs-bound}
\operatorname{tr}(\Gamma_k(C))
&=\sum_{\ell=1}^{\min(N,k)}
\sum_{\substack{J\subset\{1,\ldots,N\}\\|J|=\ell}}
\sum_{\substack{\beta\in\mathbb N^J\\
\sum_{j\in J}\beta_j=k}}
\prod_{j\in J}\lambda_j^{\beta_j}
\le
\sum_{\ell=1}^{\min(N,k)}
\sum_{\substack{J\subset\{1,\ldots,N\}\\|J|=\ell}}
\biggl(\prod_{j\in J}\lambda_j\biggr)
\sum_{\substack{\beta\in\mathbb N^J\\
\sum_{j\in J}\beta_j=k}}1\\
&=
\sum_{\ell=1}^{\min(N,k)}
\binom{k-1}{\ell-1}
\sum_{\substack{J\subset\{1,\ldots,N\}\\|J|=\ell}}
\prod_{j\in J}\lambda_j
=
\sum_{\ell=1}^{\min(N,k)}
\binom{k-1}{\ell-1}e_\ell(C).
\nonumber
\end{align}

We also have
\begin{equation}\label{eq-compression-identity}
\sum_{j=1}^m w_j e_\ell(P_{{\bf v}_j^\perp}CP_{{\bf v}_j^\perp})
=(N-\ell)e_\ell(C),
\quad \ell\in\{0,\ldots,N\}.
\end{equation}
Indeed, for $s\ge0$, the cofactor formula, applied in an orthonormal
basis whose first vector is ${\bf v}_j$, gives
\[
{\bf v}_j^*(I_N+sC)^{-1}{\bf v}_j
= \frac{\det(I_N+sP_{{\bf v}_j^\perp}CP_{{\bf v}_j^\perp})}{\det(I_N+sC)}.
\]
Hence,
\[
\det(I_N+sP_{{\bf v}_j^\perp}CP_{{\bf v}_j^\perp})
=\det(I_N+sC)\,
{\bf v}_j^*(I_N+sC)^{-1}{\bf v}_j.
\]
Summing these identities with weights $w_j$ and using the identity
$\sum_{j=1}^m w_j{\bf v}_j{\bf v}_j^*=I_N$,
we obtain
\begin{align*}
&\sum_{j=1}^m w_j\det(I_N+sP_{{\bf v}_j^\perp}CP_{{\bf v}_j^\perp})
=\det(I_N+sC)
\sum_{j=1}^m w_j
\operatorname{tr}\bigl((I_N+sC)^{-1}{\bf v}_j{\bf v}_j^*\bigr)
\\
&=
\det(I_N+sC)\operatorname{tr}\Bigl((I_N+sC)^{-1}
\sum_{j=1}^m w_j
{\bf v}_j{\bf v}_j^*\Bigr)
=
\det(I_N+sC)\operatorname{tr}(I_N+sC)^{-1}.
\end{align*}
Let $\lambda_1, \ldots, \lambda_N$ be the eigenvalues of the matrix $C$.
Then
\[
\det(I_N+sC)\operatorname{tr}(I_N+sC)^{-1}
=
\sum_{j=1}^N\prod_{r\ne j}(1+s\lambda_r)
=
\sum_{\ell=0}^N(N-\ell)e_\ell(C)s^\ell.
\]
Similarly,
\[
\det(I_N+sP_{{\bf v}_j^\perp}CP_{{\bf v}_j^\perp})
=\sum_{\ell=0}^N e_\ell(P_{{\bf v}_j^\perp}CP_{{\bf v}_j^\perp})s^\ell.
\]
Thus,
\[
\sum_{\ell=0}^N \Bigl(\sum_{j=1}^m w_je_\ell(P_{{\bf v}_j^\perp}CP_{{\bf v}_j^\perp})\Bigr)s^\ell
=
\sum_{j=1}^m w_j\det(I_N+sP_{{\bf v}_j^\perp}CP_{{\bf v}_j^\perp})
=
\sum_{\ell=0}^N(N-\ell)e_\ell(C)s^\ell.
\]
Comparing coefficients proves
\eqref{eq-compression-identity}.

Returning to $\|\Gamma_k(B_n)\|_{\mathrm{HS}}^2$,
we have $\|B_n\|_{\mathrm{op}}\le1$,
and, by Lemma \ref{lem-aux}(i),(ii) and \eqref{eq-gamma-hs-bound}, we obtain
\begin{equation}\label{eq-hs-est}
\|\Gamma_k(B_n)\|_{\mathrm{HS}}^2 =
\operatorname{tr}(\Gamma_k(B_nB_n^*))
\le
\sum_{\ell=1}^{\min(N,k)}
\binom{k-1}{\ell-1}e_\ell(C_n).
\end{equation}
Since $C_n=P_{{\bf v}_{J_n}^\perp}C_{n-1}P_{{\bf v}_{J_n}^\perp}$,
\eqref{eq-compression-identity} gives
\[
\mathbb E\bigl[e_\ell(C_n)\mid C_{n-1}\bigr]
=\left(1-\frac{\ell}{N}\right)e_\ell(C_{n-1}),
\]
and, since $C_0=I_N$, induction yields
\begin{equation}\label{eq-expect}
\mathbb E e_\ell(C_n)
=\binom N\ell\biggl(1-\frac{\ell}{N}\biggr)^n.
\end{equation}
Combining \eqref{eq-hs-est} and \eqref{eq-expect} with \eqref{eq-1-trace-est}, we obtain
\[
\operatorname{tr}e^{-tS_k}
\le
\sum_{\ell=1}^{\min(N,k)}
\binom N\ell\binom{k-1}{\ell-1}
e^{-Nt/2}\sum_{n=0}^\infty
\frac{(Nt/2)^n}{n!}\biggl(1-\frac{\ell}{N}\biggr)^n
=
\sum_{\ell=1}^{\min(N,k)}
\binom N\ell\binom{k-1}{\ell-1}e^{-t\ell/2}.
\]

Since
\[
\binom N\ell\binom{k-1}{\ell-1}
\le \frac{N^\ell k^{\ell-1}}{\ell!(\ell-1)!}\le
\frac{N(Nk)^{\ell-1}}{(\ell-1)!},
\]
we obtain
\[
\operatorname{tr}e^{-tS_k}
\le
Ne^{-t/2}\sum_{\ell=1}^\infty
\frac{(Nke^{-t/2})^{\ell-1}}{(\ell-1)!}
=Ne^{-t/2} e^{Nke^{-t/2}},
\]
which completes the proof.
\end{proof}

\section{Proof of Theorem \ref{th-uniform-signed-sum}}
\label{sect-power-sums}

Coordinates that vanish identically on $E$ can be omitted.
Let ${\bf u}_1,\ldots,{\bf u}_N$ be an orthonormal basis of $E$,
and set
\[
w_j:=\sum_{\ell=1}^N|{\bf u}_\ell(j)|^2,
\qquad
{\bf v}_j:=w_j^{-1/2}
\begin{pmatrix}
{\bf u}_1(j)\\ \vdots\\ {\bf u}_N(j)
\end{pmatrix}.
\]
Then
\[
\|{\bf v}_j\|_2=1,\quad
\sum_{j=1}^mw_j{\bf v}_j{\bf v}_j^*=I_N,
\quad \text{and}\quad
\sum_{j=1}^mw_j=N.
\]
Moreover,
\[
w_j=\sup\limits_{x\in E,\,\|x\|_2=1}|x_j|^2\le\delta.
\]
We relabel the coordinates so that $w_1\ge\ldots\ge w_m$.

For $p\in[p_0,2]$, let $w_j(p):=\delta(w_j/\delta)^{p/2}$.
These weights are nonincreasing in $j$ and we have
$w_j\le w_j(p)\le\delta$. Since the vector
$x=\sum_{\ell=1}^Ny_\ell{\bf u}_\ell$ has coordinates
$x_j=\sqrt{w_j}\langle{\bf v}_j,y\rangle$, homogeneity gives
\[
S_p(E,\zeta):=
\sup_{x\in E\cap B_p^m}
\Bigl|\sum_{j=1}^m\zeta_j|x_j|^p\Bigr|
=
\sup_{\|y\|_2=1}
\frac{
\bigl|\sum_{j=1}^m\zeta_jw_j(p)|\langle{\bf v}_j,y\rangle|^p\bigr|
}{
\sum_{j=1}^m w_j(p)|\langle{\bf v}_j,y\rangle|^p
}.
\]

For fixed $u\ge0$, let
\[
a:=\sqrt{\delta\bigl(12p_0^{-1}\log(2N)+u\bigr)},
\quad
t_k:=2\bigl(\log(2eNk(k+1))+u\bigr),
\]
and let
\[
S_k:=\sum_{j=1}^m w_jT_{k,{\bf v}_j}, \quad T_{k, {\bf v}_j} := I_{\mathcal{H}_k^N} - \Gamma_k(P_{{\bf v}_j^\perp}).
\]

Let
\[
Y_n^\pm:=
\operatorname{tr}\exp\Bigl(
-t_kS_k+2a\sum_{j=1}^n(\pm\zeta_j-a)T_{k,{\bf v}_j}
\Bigr),
\quad 0\le n\le m.
\]
Since $T_{k,{\bf v}_j}$ is an orthogonal projection,
\[
e^{sT_{k,{\bf v}_j}}
=
I_{\mathcal H_k^N}-T_{k,{\bf v}_j}
+e^sT_{k,{\bf v}_j}
\quad\forall s\in\mathbb R.
\]
Consequently, the subgaussian assumption gives
\[
\mathbb E e^{2a(\pm\zeta_j-a)}
=
e^{-2a^2}\mathbb E e^{\pm2a\zeta_j}
\le 1,
\]
and hence
\[
\mathbb E e^{2a(\pm\zeta_j-a)T_{k,{\bf v}_j}}
=
I+\bigl(\mathbb E e^{2a(\pm\zeta_j-a)}-1\bigr)T_{k,{\bf v}_j}
\preceq I
\]
and the matrix $\log \mathbb E e^{2a(\pm\zeta_j-a)T_{k,{\bf v}_j}}$
is negative
semidefinite.
By \eqref{eq-LT}
and by monotonicity of the trace exponential in the
Loewner order,
we obtain
\begin{align*}
\mathbb E[Y_n^\pm\mid \zeta_1, \ldots, \zeta_{n-1}]
&\le
\operatorname{tr}\exp\Bigl(
-t_kS_k+2a\sum_{j=1}^{n-1}(\pm\zeta_j-a)T_{k,{\bf v}_j}
+\log\mathbb E e^{2a(\pm\zeta_n-a)T_{k,{\bf v}_n}}
\Bigr)
\\
&\le
\operatorname{tr}\exp\Bigl(
-t_kS_k+2a\sum_{j=1}^{n-1}(\pm\zeta_j-a)T_{k,{\bf v}_j}
\Bigr)
=
Y_{n-1}^\pm.
\end{align*}

Thus, $(Y_n^\pm)_{n=0}^m$ is a nonnegative supermartingale.
The maximal inequality for nonnegative supermartingales
\cite[Lemma~5.18]{Durrett2012} and
Lemma~\ref{lem-heat-trace}
give
\[
\mathbb P\Bigl(\max_{0\le n\le m}Y_n^\pm>1\Bigr)
\le
\mathbb E Y_0^\pm
=
\operatorname{tr}e^{-t_kS_k}
\le Ne^{-t_k/2} e^{Nke^{-t_k/2}}\le
eNe^{-t_k/2} = \frac{e^{-u}}{2k(k+1)}.
\]
Let
\[
A:=\bigcap_{k=1}^\infty
\bigcap_{\sigma\in\{-1,1\}}
\bigcap_{n=1}^m
\biggl\{
\sum_{j=1}^n(\sigma\zeta_j-a)T_{k,{\bf v}_j}
\preceq \frac{t_k}{2a}S_k
\biggr\}.
\]
If one of these inequalities fails, the exponent in the
corresponding $Y_n^\pm$ has a positive eigenvalue, so
$Y_n^\pm>1$. The union bound over $k$ and the two signs
therefore gives
\[
\mathbb P(A^c)
\le 2\sum_{k=1}^\infty\frac{e^{-u}}{2k(k+1)}
=e^{-u}.
\]

Fix $p\in[p_0,2]$ and set $w_{m+1}(p):=0$.
On the event $A$,
\[
\sum_{j=1}^n(\pm\zeta_j-a)T_{k,{\bf v}_j}
\preceq \frac{t_k}{2a}S_k\quad \forall n\in \{1, \ldots, m\}.
\]
Multiplying these inequalities
by
$w_n(p)-w_{n+1}(p)\ge0$ and summing over $n$,
we obtain
\[
\sum_{n=1}^m\bigl(w_n(p)-w_{n+1}(p)\bigr)
\sum_{j=1}^n(\pm\zeta_j-a)T_{k,{\bf v}_j}
\preceq
\frac{t_k}{2a}S_k
\sum_{n=1}^m\bigl(w_n(p)-w_{n+1}(p)\bigr)
=
\frac{w_1(p)t_k}{2a}S_k.
\]
Reversing the order of summation on the left, we obtain
\begin{align*}
\sum_{n=1}^m\bigl(w_n(p)-w_{n+1}(p)\bigr)
\sum_{j=1}^n(\pm\zeta_j-a)T_{k,{\bf v}_j}
&=
\sum_{j=1}^m(\pm\zeta_j-a)T_{k,{\bf v}_j}
\sum_{n=j}^m\bigl(w_n(p)-w_{n+1}(p)\bigr)
\\
&=
\sum_{j=1}^m w_j(p)(\pm\zeta_j-a)T_{k,{\bf v}_j}.
\end{align*}
Since $w_1(p)\le\delta$ and $S_k\succeq0$, it follows that
\[
\pm\sum_{j=1}^m w_j(p)\zeta_jT_{k,{\bf v}_j}
\preceq
a\sum_{j=1}^m w_j(p)T_{k,{\bf v}_j}
+\frac{\delta t_k}{2a}S_k.
\]
This inequality holds simultaneously for every $p\in[p_0,2]$ on the event $A$.

For a unit vector $y\in \mathbb{R}^N$, let
\[
\pi_{p,k}(y):=
\frac{
b_{p/2,k}
\langle\psi_{k,y},S_k\psi_{k,y}\rangle_{\mathcal H_k^N}
}{
\sum_{j=1}^mw_j|\langle{\bf v}_j,y\rangle|^p
}.
\]
We have $\pi_{p,k}(y)\ge0$ and,
by representation~\eqref{eq-repres},
\begin{equation}\label{eq-probab-weights}
\sum_{k=1}^\infty\pi_{p,k}(y)=1.
\end{equation}
On the event $A$, by representations \eqref{eq-repres} and \eqref{eq-repres-random}, since $w_j(p)\ge w_j$, we then have
\[
\frac{
\bigl|\sum_{j=1}^m\zeta_jw_j(p)|\langle{\bf v}_j,y\rangle|^p\bigr|
}{
\sum_{j=1}^mw_j(p)|\langle{\bf v}_j,y\rangle|^p
}
\le
a+\frac{\delta}{2a}\sum_{k=1}^\infty\pi_{p,k}(y)t_k.
\]
It remains to estimate the last sum.

We have
\begin{equation}\label{eq-pi-1}
\sum_{j=1}^mw_j|\langle{\bf v}_j,y\rangle|^p
\ge
\sum_{j=1}^mw_j|\langle{\bf v}_j,y\rangle|^2
=1,
\end{equation}
and, by
Lemma \ref{lem-aux}(iv),
\begin{equation}\label{eq-pi-2}
\langle\psi_{k,y},S_k\psi_{k,y}\rangle_{\mathcal H_k^N}
=
\sum_{j=1}^m w_j\bigl[1-\bigl(1-\langle y, {\bf v}_j\rangle^2\bigr)^k\bigr]
\le\sum_{j=1}^m w_j=N.
\end{equation}
For the binomial coefficients, for every $\alpha\in (0, 1]$, we have
\[
b_{\alpha, k}
=\frac{\alpha}{k}
\prod_{\ell=1}^{k-1}\Bigl(1-\frac{\alpha}{\ell}\Bigr)
=
\prod_{\ell=1}^{k-1}\Bigl(1-\frac{\alpha}{\ell}\Bigr) -
\prod_{\ell=1}^{k}\Bigl(1-\frac{\alpha}{\ell}\Bigr).
\]
Therefore, using this representation and  \eqref{eq-coeff-equality}, we conclude
\begin{equation}\label{eq-pi-3}
\sum_{k=n+1}^\infty b_{\alpha, k}
=
1- \sum_{k=1}^nb_{\alpha, k}
=
\prod_{\ell=1}^n\left(1-\frac{\alpha}{\ell}\right)
\le
\exp\biggl(-\alpha\sum_{\ell=1}^n\frac1\ell\biggr)
\le(n+1)^{-\alpha}\quad \forall n\in\mathbb{N}\cup\{0\}.
\end{equation}
Consequently, by \eqref{eq-probab-weights}, \eqref{eq-pi-1}, \eqref{eq-pi-2}, and \eqref{eq-pi-3},
\[
\sum_{k=n+1}^\infty\pi_{p,k}(y)
\le
\min\{1,N(n+1)^{-p/2}\}\quad \forall n\in \mathbb{N}\cup\{0\}.
\]
Hence,
\begin{align*}
\sum_{k=1}^\infty\pi_{p,k}(y)\log k
&=\int_0^\infty
\sum_{\log k>t}\pi_{p,k}(y)\,dt
=
\int_0^\infty
\sum_{k=\lfloor e^t\rfloor+1}^\infty\pi_{p,k}(y)\,dt
\\
&\le\int_0^\infty\min\{1,N(\lfloor e^t\rfloor+1)^{-p/2}\}\,dt
\le\int_0^\infty\min\{1,Ne^{-  pt/2}\}\,dt
=\frac{2}{p}(\log N+1).
\end{align*}
Since $k(k+1)\le 2k^2$ and $\sum_{k=1}^\infty\pi_{p,k}(y)=1$,
the definition of $t_k$ gives
\begin{align*}
\frac12\sum_{k=1}^\infty\pi_{p,k}(y)t_k
&\le
\log(4eN)+2\sum_{k=1}^\infty\pi_{p,k}(y)\log k+u
\\
&\le
\log(4eN)+\frac4p(\log N+1)+u
\le
12p_0^{-1}\log(2N)+u.
\end{align*}
In the last step, we used
$\log(4eN)\le4\log(2N)$ and $\log N+1\le2\log(2N)$.
Therefore, on the event $A$,
\[
W:=\sup_{p_0\le p\le2}S_p(E,\zeta)
\le
a+\frac{\delta}{a}\bigl(12p_0^{-1}\log(2N)+u\bigr)
= 2\sqrt{\delta\bigl(12p_0^{-1}\log(2N)+u\bigr)}.
\]

It follows that the random variable
$U:=\bigl(W^2/(4\delta)-12p_0^{-1}\log(2N)\bigr)_+$
satisfies
\[
\mathbb P(U>u)\le e^{-u} \quad  \forall u\ge0.
\]
For every $n\in \mathbb{N}$,
we obtain
\[
\mathbb E U^n
=n\int_0^\infty u^{n-1}\mathbb P(U>u)\,du
\le n!.
\]
Therefore,
\[
\bigl(\mathbb{E}W^{2n}\bigr)^{1/n}
\le 4\delta\bigl(12p_0^{-1}\log(2N)+n\bigr)
\le 64\delta\bigl(p_0^{-1}\log(2N)+n\bigr).
\]
For any $r\ge1$, taking $n=\lceil r/2\rceil$,
since $r\le2n$ and $n\le r$, we obtain
\[
\bigl(\mathbb E W^r\bigr)^{1/r}
\le \bigl(\mathbb{E}W^{2n}\bigr)^{1/(2n)}
\le8\sqrt{\delta\bigl(p_0^{-1}\log(2N)+r\bigr)},
\]
as announced.
\qed

\section{Proof of Theorem~\ref{th-uniform-sampling}}
\label{sect-sampling}

We first consider the real case and assume that
\[
m\ge 2^{19}p_0^{-1}\varepsilon^{-2}KN\log(2N).
\]

Let $\eta_1,\ldots,\eta_m$ be independent Rademacher variables,
independent of $\xi:=(\xi^1,\ldots,\xi^m)$, and let
\[
Z_\xi:=
\sup_{p_0\le p\le2}
\sup_{\substack{f\in F_N\\\|f\|_{L_p(\mu)}=1}}
\biggl|\frac1m\sum_{j=1}^m|f(\xi^j)|^p-1\biggr|
\quad\text{and}\quad
R:=
\frac1m\sup_{p_0\le p\le2}
\sup_{\substack{f\in F_N\\\|f\|_{L_p(\mu)}=1}}
\biggl|\sum_{j=1}^m\eta_j|f(\xi^j)|^p\biggr|.
\]
The standard symmetrization argument gives
\begin{equation}\label{eq-symm}
\bigl(\mathbb E Z_\xi^r\bigr)^{1/r} \le2\bigl(\mathbb E R^r\bigr)^{1/r}
\quad\forall r\ge1.
\end{equation}

Let $\widetilde\xi$ be an independent copy of $\xi$, also
independent of the signs, and let
\[
A:=
\biggl\{
\frac1m\sum_{j=1}^m|f(\widetilde\xi^j)|^2
\ge \frac12\|f\|_{L_2(\mu)}^2
\quad\text{for every }f\in F_N
\biggr\}.
\]
By
\cite[Corollary~5.2]{Tropp2012}
(see also \cite[Remark 3.5 and the proof of Theorem 3.2]{DKT26}),
\begin{equation}\label{eq-probab-est}
\mathbb P(A)
\ge1-N\exp\Bigl(-\frac{m}{8KN}\Bigr)\ge \frac{1}{2}.
\end{equation}
For fixed samples $\widetilde{\xi}\in A$ and $\xi$, consider the space
\[
E:=
\bigl\{
(f(\xi^1),\ldots,f(\xi^m),
f(\widetilde\xi^1),\ldots,f(\widetilde\xi^m))\colon
f\in F_N
\bigr\}\subset\mathbb R^{2m}.
\]
The evaluation map is injective on $A$, so $\dim E=N$.
For $x\in E$ and for the corresponding function $f\in F_N$, we have
\[
\|x\|_\infty^2
\le KN\|f\|_{L_2(\mu)}^2
\le\frac{2KN}{m}\sum_{j=1}^m|f(\widetilde\xi^j)|^2
\le\frac{2KN}{m}\|x\|_2^2.
\]
Applying Theorem~\ref{th-uniform-signed-sum} with the
subgaussian coefficients
$(\eta_1,\ldots,\eta_m,0,\ldots,0)$,
we obtain
\begin{equation}\label{eq-r-moment-1}
\biggl[
\mathbb E_\eta
\biggl(
\sup_{p_0\le p\le2}\sup_{x\in E\cap B_p^{2m}}
\Bigl|\sum_{j=1}^m\eta_j|x_j|^p\Bigr|
\biggr)^r
\biggr]^{1/r}
\le
8\sqrt{\frac{2KN}{m}\bigl(p_0^{-1}\log(2N)+r\bigr)}.
\end{equation}

For each $p\in[p_0,2]$ and each $f\in F_N$ with
$\|f\|_{L_p(\mu)}=1$,
\[
\sum_{j=1}^m
\bigl(|f(\xi^j)|^p+|f(\widetilde\xi^j)|^p\bigr)
\le m(2+Z_\xi+ Z_{\widetilde{\xi}}).
\]
Therefore, on $A$,
\[
R\le
(2+Z_\xi+ Z_{\widetilde{\xi}})
\sup_{p_0\le p\le 2}\sup_{x\in E\cap B_p^{2m}}
\biggl|\sum_{j=1}^m\eta_j|x_j|^p\biggr|,
\]
and, by \eqref{eq-r-moment-1},
\[
\bigl[\mathbb E_\eta R^r\bigr]^{1/r}
\le 8(2+Z_\xi+ Z_{\widetilde{\xi}})
\sqrt{\frac{2KN}{m}\bigl(p_0^{-1}\log(2N)+r\bigr)}
\quad\text{for } \widetilde{\xi}\in A.
\]
Since $\mathbf1_A$ is a function of variables $\widetilde{\xi}$,
from \eqref{eq-probab-est},
we obtain
\[
\bigl(\mathbb{E}_{\widetilde{\xi}, \xi, \eta}\bigl[\mathbf1_AR^r\bigr]\bigr)^{1/r}
= \mathbb{P}(A)^{1/r}
\bigl(\mathbb{E}_{\xi, \eta}R^r\bigr)^{1/r}
\ge \frac{1}{2}\bigl(\mathbb{E}_{\xi, \eta}R^r\bigr)^{1/r}.
\]
Hence, by \eqref{eq-symm},
\begin{align*}
\bigl(\mathbb E Z_\xi^r\bigr)^{1/r} &\le
2\bigl(\mathbb E_{\xi, \eta} R^r\bigr)^{1/r}
\le 4 \bigl(\mathbb{E}_{\widetilde{\xi}, \xi}\mathbb{E}_\eta\bigl[\mathbf1_AR^r\bigr]\bigr)^{1/r}
\\
&\le
32\sqrt{\frac{2KN}{m}\bigl(p_0^{-1}\log(2N)+r\bigr)}
\bigl(\mathbb{E}_{\widetilde{\xi}, \xi}(2+Z_\xi+ Z_{\widetilde{\xi}})^r\bigr)^{1/r}
\\
&\le
64
\sqrt{\frac{2KN}{m}\bigl(p_0^{-1}\log(2N)+r\bigr)}
\bigl(1+(\mathbb{E}Z_\xi^r)^{1/r}
\bigr).
\end{align*}

Let
\[
r= 2^{-18}\frac{m\varepsilon^2}{KN}.
\]
Since
\[
m\ge 2^{19}p_0^{-1}\varepsilon^{-2}KN\log(2N),
\]
we have $r\ge 2p_0^{-1}\log(2N)\ge 1$ and
\[
64
\sqrt{\frac{2KN}{m}\bigl(p_0^{-1}\log(2N)+r\bigr)}
=64\sqrt{\frac{2p_0^{-1}KN\log(2N)}m+2^{-17}\varepsilon^2}
\le
64\sqrt{2^{-16}\varepsilon^2}
=\frac{\varepsilon}{4}.
\]
Therefore,
\[
\bigl(\mathbb E Z_\xi^r\bigr)^{1/r}\le \frac{\varepsilon}{4}
\bigl(1+(\mathbb{E}Z_\xi^r)^{1/r}
\bigr)
\le \frac{\varepsilon}{4}+ \frac{1}{2}(\mathbb{E}Z_\xi^r)^{1/r}.
\]
Thus,
\[
\bigl(\mathbb E Z_\xi^r\bigr)^{1/r}\le \frac{\varepsilon}{2},
\]
and, by Markov's inequality,
\[
\mathbb P(Z_\xi>\varepsilon)
\le2^{-r}
\le\exp\Bigl(-2^{-19}\frac{m\varepsilon^2}{KN}\Bigr).
\]
This proves the simultaneous estimate in the real case.

We now deduce the complex case from the real case.
Let $F_N$ be a complex $N$-dimensional space satisfying
\[
\|f\|_\infty\le\sqrt {KN}\,\|f\|_{L_2(\mu)}
\quad\forall f\in F_N.
\]
Let $\nu$ be the uniform probability measure on $[0,2\pi]$.
We consider the real linear space
\[
Y:=\{g_f\colon f\in F_N\},\quad g_f(\omega,\theta):=
\operatorname{Re}\bigl(e^{i\theta}f(\omega)\bigr),
\]
of dimension $2N$
on the probability space
$(\Omega\times[0,2\pi],\mu\otimes\nu)$.
We note that
\[
\|g_f\|_\infty=\|f\|_\infty
\quad\text{and}\quad
\|g_f\|_{L_2(\mu\otimes\nu)}^2
=\frac12\|f\|_{L_2(\mu)}^2.
\]
Consequently,
\[
\|g_f\|_\infty
\le\sqrt{2KN}\,\|g_f\|_{L_2(\mu\otimes\nu)}.
\]
Assume that
\[
m\ge2^{21}p_0^{-1}\varepsilon^{-2}KN\log(2N).
\]
Since $\log(4N)\le2\log(2N)$, we have
\[
m\ge2^{19}p_0^{-1}\varepsilon^{-2}(2KN)\log(4N).
\]
Applying the real case to $Y$ with $2N$ in place of $N$, for independent samples
$(\xi^j,\theta_j)$ with distribution $\mu\otimes\nu$,
with probability at least
$1-\exp\bigl(-2^{-20}\frac{m\varepsilon^2}{KN}\bigr)$,
we have
\[
(1-\varepsilon)\|g_f\|_{L_p(\mu\otimes\nu)}^p
\le
\frac1m\sum_{j=1}^m
\bigl|\operatorname{Re}
(e^{i\theta_j}f(\xi^j))\bigr|^p
\le
(1+\varepsilon)\|g_f\|_{L_p(\mu\otimes\nu)}^p
\quad\forall f\in F_N\quad \forall p\in[p_0, 2].
\]
In particular, for every $t\in[0,2\pi]$,
\[
(1-\varepsilon)\|g_{e^{it}f}\|_{L_p(\mu\otimes\nu)}^p
\le
\frac1m\sum_{j=1}^m
\bigl|\operatorname{Re}
(e^{i(\theta_j+t)}f(\xi^j))\bigr|^p
\le
(1+\varepsilon)\|g_{e^{it}f}\|_{L_p(\mu\otimes\nu)}^p
\quad\forall f\in F_N\quad \forall p\in[p_0, 2].
\]
Integrating with respect to $t$ over $[0,2\pi]$,
since
\[
\int_0^{2\pi}
\bigl|\operatorname{Re}(e^{i(t+\theta)}w)\bigr|^p\,dt
=
|w|^p\int_0^{2\pi}|\cos t|^p\,dt
\quad\forall w\in\mathbb C\quad
\forall\theta\in[0,2\pi],
\]
we obtain
\[
(1-\varepsilon)\|f\|_{L_p(\mu)}^p
\le
\frac1m\sum_{j=1}^m|f(\xi^j)|^p
\le
(1+\varepsilon)\|f\|_{L_p(\mu)}^p
\quad\forall f\in F_N\quad \forall p\in[p_0, 2]
\]
with probability at least
$1-\exp\bigl(-2^{-20}\frac{m\varepsilon^2}{KN}\bigr)$,
as claimed.
\qed

\section{Weighted discretization and embeddings of subspaces}
\label{sect-embedding}

A standard application of Lewis's change of density theorem
\cite{Lewis1978}, extended to $0<p<1$ by Schechtman and
Zvavitch \cite[Theorem~2.1]{SchechtmanZvavitch2001},
yields the following weighted discretization result.

\begin{corollary}\label{cor-weighted-discr}
Let $p\in(0,2]$, and let $F_N$ be an $N$-dimensional subspace of $L_p(\mu)$.
For every $\varepsilon\in(0,1)$ and every $m\in\mathbb N$
satisfying
\[
m\ge
2^{21}\max\{1,p^{-1}\}\varepsilon^{-2}N\log(2N),
\]
there exist points $\xi^1,\ldots,\xi^m\in\Omega$ and weights
$\lambda_1,\ldots,\lambda_m\ge0$ such that
\[
(1-\varepsilon)\|f\|_{L_p(\mu)}^p
\le
\sum_{j=1}^m\lambda_j|f(\xi^j)|^p
\le
(1+\varepsilon)\|f\|_{L_p(\mu)}^p
\quad\forall f\in F_N.
\]
Here point evaluations are understood using fixed measurable
representatives of a basis of $F_N$, extended to $F_N$ by linearity.
\end{corollary}

The corresponding embedding result follows by adjusting
the accuracy parameter.

\begin{corollary}\label{cor-Lp-embedding}
Let $p\in(0,2]$ and $\varepsilon\in(0,1)$.
Every $N$-dimensional subspace $E$ of $L_p(\mu)$
admits a linear embedding
$T\colon E\longrightarrow\ell_p^m$
such that
\[
\|T\|\,\|T^{-1}\|\le1+\varepsilon
\]
and
\[
m\le Cp^{-3}\varepsilon^{-2}N\log(2N).
\]
Here $C>0$ is an absolute constant and $T^{-1}$ is defined
on $T(E)$.
\end{corollary}

\begin{proof}
Put
\[
a:=\min\{p,1\},
\quad
\eta:=\frac{a\varepsilon}{8},
\quad\text{and}\quad
m:=\bigl\lceil
2^{27}a^{-3}\varepsilon^{-2}N\log(2N)
\bigr\rceil.
\]
Applying Corollary~\ref{cor-weighted-discr} with $\eta$
in place of $\varepsilon$, we obtain points
$\xi^1,\ldots,\xi^m$ and nonnegative weights $\lambda_1,\ldots,\lambda_m$
such that
\[
(1-\eta)\|f\|_{L_p(\mu)}^p
\le
\sum_{j=1}^m\lambda_j|f(\xi^j)|^p
\le
(1+\eta)\|f\|_{L_p(\mu)}^p
\quad\forall f\in E.
\]
Therefore, the operator
\[
Tf:=
\bigl(\lambda_1^{1/p}f(\xi^1),\ldots,
\lambda_m^{1/p}f(\xi^m)\bigr)
\]
is injective and
\[
\|T\|\,\|T^{-1}\|
\le
\biggl(\frac{1+\eta}{1-\eta}\biggr)^{1/p}.
\]
Since $\eta<1/8<1/2$ and $\log(1+\varepsilon)\ge\varepsilon/2$,
we have
\[
\frac1p\log\frac{1+\eta}{1-\eta}
\le
\frac{4\eta}{p}
=
\frac{a\varepsilon}{2p}
\le
\frac{\varepsilon}{2}
\le
\log(1+\varepsilon).
\]
Consequently, $\|T\|\,\|T^{-1}\|\le1+\varepsilon$.
Finally, since
$a^{-3}\le8p^{-3}$
for $0<p\le2$,
the assertion follows.
\end{proof}

\section*{Use of AI Tools}

The authors acknowledge the use of ChatGPT 6 Astra during the preparation of the manuscript to assist with the development of the main mathematical argument, and of ChatGPT 6.1 Sol for language and stylistic editing. 
All AI-generated material used was reviewed and revised by the authors.
The authors take full responsibility for the content of the paper.


\begin{thebibliography}{99}
\normalsize

\bibitem{BSS2012}
J.~Batson, D.~A.~Spielman, and N.~Srivastava,
Twice-Ramanujan sparsifiers,
\emph{SIAM J. Comput.}
\textbf{41} (2012), no.~6, 1704--1721.

\bibitem{BLM1989}
J.~Bourgain, J.~Lindenstrauss, and V.~Milman,
Approximation of zonoids by zonotopes,
\emph{Acta Math.}
\textbf{162} (1989), 73--141.

\bibitem{CohenPeng2015}
M.~B.~Cohen and R.~Peng,
$\ell_p$ row sampling by Lewis weights,
in \emph{Proceedings of the Forty-Seventh Annual ACM Symposium
on Theory of Computing (STOC 2015)},
ACM, New York, 2015, pp.~183--192.
Full version: arXiv:1412.0588.

\bibitem{DKM}
F.~Dai, E.~D.~Kosov, and N.~Murasko,
Empirical approximation of $L_p$ norms for $p\in(2, \infty)$,
preprint, arXiv:2611.08896v1  (2026).

\bibitem{DKT}
F.~Dai, E.~D.~Kosov, and V.~N.~Temlyakov,
Some improved bounds in sampling discretization of integral norms,
\emph{J. Funct. Anal.}
\textbf{285} (2023), no.~4, article no.~109951.

\bibitem{DKT26}
F.~Dai, E.~D.~Kosov, and V.~N.~Temlyakov,
A survey of sampling discretization of integral and uniform norms,
\emph{Surv. Approx. Theory}
\textbf{10} (2026), 1--42.

\bibitem{Durrett2012}
R.~Durrett,
\emph{Essentials of Stochastic Processes},
2nd ed., Springer Texts in Statistics, Springer, New York, 2012.

\bibitem{JLLS2023}
A.~Jambulapati, J.~R.~Lee, Y.~P.~Liu, and A.~Sidford,
Sparsifying sums of norms,
in \emph{Proceedings of the 64th IEEE Annual Symposium
on Foundations of Computer Science (FOCS 2023)},
IEEE, 2023, pp.~1953--1962.
Full version: arXiv:2305.09049.

\bibitem{JS82}
W.~B.~Johnson and G.~Schechtman,
Embedding $\ell_p^m$ into $\ell_1^n$,
\emph{Acta Math.}
\textbf{149} (1982), 71--85.

\bibitem{JS}
W.~B.~Johnson and G.~Schechtman,
Finite dimensional subspaces of $L_p$,
in \emph{Handbook of the Geometry of Banach Spaces},
W.~B.~Johnson and J.~Lindenstrauss (eds.),
vol.~1, North-Holland, Amsterdam, 2001, pp.~837--870.

\bibitem{KKLT}
B.~S.~Kashin, E.~D.~Kosov, I.~Limonova, and V.~N.~Temlyakov,
Sampling discretization and related problems,
\emph{J. Complexity}
\textbf{71} (2022), article no.~101653.

\bibitem{Lewis1978}
D.~R.~Lewis,
Finite dimensional subspaces of $L_p$,
\emph{Studia Math.}
\textbf{63} (1978), no.~2, 207--212.

\bibitem{ManojOvsiankin2025}
N.~S.~Manoj and M.~Ovsiankin,
The change-of-measure method, block Lewis weights,
and approximating matrix block norms,
in \emph{Proceedings of the 2025 Annual ACM-SIAM Symposium
on Discrete Algorithms (SODA 2025)},
SIAM, 2025, pp.~2252--2300.
Full version: arXiv:2311.10013.

\bibitem{ReisRothvoss2026}
V.~Reis and T.~Rothvoss,
Linear-size $\ell_1$ sparsifiers,
preprint, arXiv:2606.28147v2 (2026).

\bibitem{Rud}
M.~Rudelson,
Almost orthogonal submatrices of an orthogonal matrix,
\emph{Israel J. Math.}
\textbf{111} (1999), 143--155.

\bibitem{Rud99}
M.~Rudelson,
Random vectors in the isotropic position,
\emph{J. Funct. Anal.}
\textbf{164} (1999), no.~1, 60--72.

\bibitem{Schechtman1987}
G.~Schechtman,
More on embedding subspaces of $L_p$ in $\ell_r^n$,
\emph{Compos. Math.}
\textbf{61} (1987), no.~2, 159--169.

\bibitem{Schechtman2011}
G.~Schechtman,
Tight embedding of subspaces of $L_p$ in $\ell_p^n$ for even $p$,
\emph{Proc. Amer. Math. Soc.}
\textbf{139} (2011), no.~12, 4419--4421.

\bibitem{SchechtmanZvavitch2001}
G.~Schechtman and A.~Zvavitch,
Embedding subspaces of $L_p$ into $\ell_p^N$, $0<p<1$,
\emph{Math. Nachr.}
\textbf{227} (2001), no.~1, 133--142.

\bibitem{Tal90}
M.~Talagrand,
Embedding subspaces of $L_1$ into $\ell_1^N$,
\emph{Proc. Amer. Math. Soc.}
\textbf{108} (1990), no.~2, 363--369.

\bibitem{Talagrand1995}
M.~Talagrand,
Embedding subspaces of $L_p$ in $\ell_p^N$,
in \emph{Geometric Aspects of Functional Analysis
(Israel, 1992--1994)},
J.~Lindenstrauss and V.~Milman (eds.),
Oper. Theory Adv. Appl., vol.~77,
Birkh\"auser, Basel, 1995, pp.~311--326.

\bibitem{Talagrand2014}
M.~Talagrand,
\emph{Upper and Lower Bounds for Stochastic Processes:
Modern Methods and Classical Problems},
Springer, Berlin--Heidelberg, 2014.

\bibitem{Tropp2012}
J.~A.~Tropp,
User-friendly tail bounds for sums of random matrices,
\emph{Found. Comput. Math.}
\textbf{12} (2012), no.~4, 389--434.

\bibitem{Tropp2015}
J.~A.~Tropp,
An introduction to matrix concentration inequalities,
\emph{Found. Trends Mach. Learn.}
\textbf{8} (2015), nos.~1--2, 1--230.

\bibitem{Zva07}
A.~Zvavitch, More on embedding subspaces of $L_p$ into $\ell_p^N$,
$0<p<1$, in: \emph{Geometric Aspects of Functional Analysis:
Israel Seminar 1996--2000}, Springer, Berlin--Heidelberg,
2007, pp.~269--280.

\end{thebibliography}
\end{document}